\documentclass[12pt]{article}
\usepackage{amssymb,amsmath}
\usepackage{hyperref}
\usepackage{graphicx}
\usepackage{color}
\usepackage{chngcntr}

\begin{document}

\title{\LARGE\bf A new application methodology \\
of the Fourier transform for rational \\
approximation of the complex \\
error function}

%\bigskip
\author{
\normalsize\bf S. M. Abrarov\footnote{\scriptsize{Dept. Earth and Space Science and Engineering, York University, Toronto, Canada, M3J 1P3.}}\, and B. M. Quine$^{*}$\footnote{\scriptsize{Dept. Physics and Astronomy, York University, Toronto, Canada, M3J 1P3.}}}

\date{January 31, 2016}
\maketitle
%\vspace{1cm}%\bigskip

\begin{abstract}
This paper presents a new approach in application of the Fourier transform to the complex error function resulting in an efficient rational approximation. Specifically, the computational test shows that with only $17$ summation terms the obtained rational approximation of the complex error function provides the average accuracy ${10^{ - 15}}$ over the most domain of practical importance $0 \le x \le 40,000$ and ${10^{ - 4}} \le y \le {10^2}$ required for the HITRAN-based spectroscopic applications. Since the rational approximation does not contain trigonometric or exponential functions dependent upon the input parameters $x$ and $y$, it is rapid in computation. Such an example demonstrates that the considered methodology of the Fourier transform may be advantageous in practical applications.
\vspace{0.25cm}
\\
\noindent {\bf Keywords:} Fourier transform, complex error function, Faddeeva function, rational approximation
\vspace{0.25cm}
\end{abstract}

\section{Introduction}
The forward and inverse Fourier transforms can be defined as \cite{Bracewell2000, Hansen2014}
\[
\label{eq_1a}\tag{1a}
F\left( \nu  \right) = {\cal{F}}\left\{ {f\left( t \right)} \right\}\left( \nu  \right) = \int\limits_{ - \infty }^\infty  {f\left( t \right){e^{ - 2\pi i \nu t}}dt}
\]
and
\[
\label{eq_1b}\tag{1b}
f\left( t \right) = {\cal{F}}^{-1}\left\{ {F\left( \nu  \right)} \right\}\left( t \right) = \int\limits_{ - \infty }^\infty  {F\left( \nu  \right){e^{2\pi i \nu t}}d\nu },
\]
respectfully. Approximation theory based on the Fourier trigonometric series for functions or signals remains a topical subject in mathematical analysis and many new efficient methodologies have been reported in the recent scientific literature (see for example \cite{Mishra2012, Mishra2013, Boyd2013}).

In our recent publication \cite{Abrarov2015a} we have shown that a sampling with the Gaussian function of the kind $h{e^{ - {{\left( {t/c} \right)}^2}}}/ \left( {c\sqrt \pi  } \right)$ leads to the trigonometric approximations for the forward
\[
\label{eq_2a}\tag{2a}
F\left( \nu  \right) = {\cal{F}}\left\{ {f\left( t \right)} \right\}\left( \nu  \right) \approx h{e^{ - {{\left( {\pi c\nu } \right)}^2}}}\sum\limits_{n =  - N}^N {f\left( {nh} \right){e^{ - 2\pi i\nu nh}}}
\]
and inverse Fourier transforms
\[
\label{eq_2b}\tag{2b}
f\left( t \right) = {\cal{F}}^{-1}\left\{ {F\left( \nu  \right)} \right\}\left( t \right) \approx h{e^{ - {{\left( {\pi ct} \right)}^2}}}\sum\limits_{n =  - N}^N {F\left( {nh} \right){e^{2\pi itnh}}},
\]
where $h$ is the step between two adjacent sampling points and $c$ is the fitting parameter, $e^{ - 2\pi i\nu nh} = \cos \left( 2\pi i\nu nh \right) - i\sin \left( 2\pi i\nu nh \right)$ and $e^{ 2\pi itnh} = \cos \left( 2\pi itnh \right) + i\sin \left( 2\pi itnh \right)$. The parameters $h$, $c$ and $N$ in the equations \eqref{eq_2a} and \eqref{eq_2b} may be the same in the forward and inverse Fourier transforms only when we imply the most favorable conditions $h <  < 1$, $c <  < 1$ and $N >  > 1$. In practical tasks, however, these conditions may be compromised in order to reduce the number of the summation terms. As a result, these parameters may not be necessarily equal to each other in the forward and inverse Fourier transforms. Consequently, it is convenient to rewrite two equations above in form
$$
F\left( \nu  \right) = {\cal{F}}\left\{ {f\left( t \right)} \right\}\left( \nu  \right) \approx {h_f}{e^{ - {{\left( {\pi {c_f}\nu } \right)}^2}}}\sum\limits_{m =  - M}^M {f\left( {m{h_f}} \right){e^{ - 2\pi i\nu m{h_f}}}}
$$
and
$$
f\left( t \right) = {\cal{F}}^{-1}\left\{ {F\left( \nu  \right)} \right\}\left( t \right) \approx {h_i}{e^{ - {{\left( {\pi {c_i}t} \right)}^2}}}\sum\limits_{n =  - N}^N {F\left( {n{h_i}} \right){e^{2\pi itnh{_i}}}},
$$
where ${h_f}$, ${c_f}$, $M$ and ${h_i}$, ${c_i}$, $N$ are the steps, the fitting parameters and the integers corresponding to the forward and inverse Fourier transforms, respectively. 

The presence of the damping functions ${e^{ - {{\left( {\pi {c_f}\nu } \right)}^2}}}$ and ${e^{ - {{\left( {\pi {c_i}t} \right)}^2}}}$ in the equations above excludes periodicity of the approximated functions $f\left( t \right)$ and $F\left( \nu  \right)$. Consequently, a solitary wavelet (or non-periodic pulse) can be effectively approximated in the Fourier transform. However, when we take ${c_f} = {c_i} = 0$, the right side of these equations become periodic with corresponding periods $1/{h_f}$, $1/{h_i}$ and represent the Fourier-type expansion series as follows
\[
\label{eq_3a}\tag{3a}
F\left( \nu  \right) \approx {h_f}\sum\limits_{m =  - M}^M {f\left( {m{h_f}} \right){e^{ - 2\pi i\nu m{h_f}}}},	\qquad - \frac{1}{{2{h_f}}} \le \nu  \le \frac{1}{{2{h_f}}},
\]
and
\[
\label{eq_3b}\tag{3b}
f\left( t \right) \approx {h_i}\sum\limits_{n =  - N}^N {F\left( {n{h_i}} \right){e^{2\pi itn{h_i}}}}, \qquad - \frac{1}{{2{h_i}}} \le t \le \frac{1}{{2{h_i}}}.
\]
It should be noted that if the integral \eqref{eq_1a} is not analytically integrable, then the function $f\left( t \right)$ can be approximated by substituting equation \eqref{eq_3a} into \eqref{eq_3b}. This substitution yields
\small
\[
\begin{aligned}
f\left( t \right) \approx {h_i}\sum\limits_{n =  - N}^N {\underbrace {\left[ {{h_f}\sum\limits_{m =  - M}^M {f\left( {m{h_f}} \right){e^{ - 2\pi in{h_i}m{h_f}}}} } \right]}_{ \approx F\left( {n{h_i}} \right)}{e^{2\pi itn{h_i}}}} \\
 = {h_i}{h_f}\sum\limits_{n =  - N}^N {\sum\limits_{m =  - M}^M {f\left( {m{h_f}} \right){e^{2\pi in{h_i}\left( {t - m{h_f}} \right)}}} } , &  &  - \frac{1}{{2{h_i}}} \le t \le \frac{1}{{2{h_i}}}.
\end{aligned}
\]
\normalsize

In this work we show a new application methodology of the Fourier transform to the complex error function. Due to representation of the complex error function as a rational approximation, it is rapid in computation. Furthermore, with only 17 summation terms the obtained rational approximation of the complex error function provides accuracy ${10^{ - 15}}$ over the most domain of practical importance $0 \le x \le 40,000 \cap {10^{ - 4}} \le y \le {10^2}$ required for applications utilizing the HITRAN molecular spectroscopic database \cite{Rothman2013}.

\section{Derivation}

\subsection{Function overview}
The complex error function, also known as the Faddeeva function or the Kramp function, can be defined as \cite{Faddeyeva1961, Gautschi1970, Abramowitz1972, Armstrong1972, Schreier1992}
$$
w\left( z \right) = {e^{ - {z^2}}}\left( {1 + \frac{{2i}}{{\sqrt \pi  }}\int\limits_0^z {{e^{{t^2}}}dt} } \right).
$$
where $z = x + iy$ is the complex argument. The complex error function is a solution of the differential equation \cite{Schreier1992}
$$
w'\left( z \right) + 2zw\left( z \right) = \frac{{2i}}{{\sqrt \pi  }},
$$
with initial condition $w\left( 0 \right) = 1.$

The complex error function is closely related to a family of the special functions. Among them the most important one is the complex probability function \cite{Armstrong1972, Schreier1992, Weideman1994}
$$
W\left( z \right) = PV\frac{i}{\pi }\int\limits_{ - \infty }^\infty  {\frac{{{e^{ - {t^2}}}}}{{z - t}}dt}
$$
or
$$
W\left( {x,y} \right) = PV\frac{i}{\pi }\int\limits_{ - \infty }^\infty  {\frac{{{e^{ - {t^2}}}}}{{\left( {x + iy} \right) - t}}dt}.
$$
This principal value integral implies that the complex probability function has no discontinuity at $y = 0$ and $x = t$. In particular,
\setcounter{equation}{3}
\begin{equation}\label{eq_4}
\lim W\left( {x,y \to 0} \right) = {e^{ - {x^2}}} + \frac{{2i}}{{\sqrt \pi  }}{\rm{daw}}\left( x \right),
\end{equation}
where ${\rm{daw}}\left( x \right)$ is the Dawson\text{'}s integral that will be briefly introduced later. There is a direct relationship between complex error function and complex probability function \cite{Armstrong1972, Schreier1992}
\begin{equation}\label{eq_5}
W\left( z \right) = w\left( z \right), \qquad {\mathop{\rm Im}\nolimits} \left[ z \right] \ge 0.
\end{equation}

The real part of the complex probability function, denoted as $K\left( {x,y} \right)$,  is known as the Voigt function. Mathematically, the Voigt function represents a convolution integral of the Gaussian and Lorentzian distributions \cite{Armstrong1972, Schreier1992, Letchworth2007, Pagnini2010, Abrarov2010}
$$
K\left( {x,y} \right) = PV\frac{y}{\pi }\int\limits_{ - \infty }^\infty  {\frac{{{e^{ - {t^2}}}}}{{{y^2} + {{\left( {x - t} \right)}^2}}}dt,}
$$
where the principal value integral also implies that it has no discontinuity at $y = 0$ and $x = t$. Specifically, from equation \eqref{eq_4} it follows that
$$
\lim K\left( {x,y \to 0} \right) = {e^{ - {x^2}}}.
$$
At non-negative argument $y$ the real part of the complex error function is also the Voigt function in accordance with identity \eqref{eq_5}. The Voigt function is widely used in many spectroscopic applications as it describes the line broadening effects \cite{Edwards1992, Quine2002, Christensen2012, Berk2013, Quine2013}.  Therefore, the application of the complex error function is very significant in quantitative spectroscopy.

Other closely related functions are the error function of complex argument \cite{Schreier1992}
$$
w\left( z \right) = {e^{ - {z^2}}}{\rm{erfc}}\left( { - iz} \right) = {e^{ - {z^2}}}\left[ {1 - {\rm{erf}}\left( { - iz} \right)} \right] \quad  \Leftrightarrow  \quad {\rm{erf}}\left( z \right) = 1 - {e^{ - {z^2}}} w\left( {iz} \right),
$$
the plasma dispersion function \cite{Fried1961}
$$
Z\left( z \right) = PV\frac{1}{{\sqrt \pi  }}\int\limits_{ - \infty }^\infty  {\frac{{{e^{ - {t^2}}}}}{{t - z}}dt = i\sqrt \pi  w\left( z \right)}
$$
the Dawson\text{'}s integral \cite{Cody1970, McCabe1974, Rybicki1989, Boyd2008, Abrarov2015b}
$$
{\rm{daw}}\left( z \right) = {e^{ - {z^2}}}\int\limits_0^z {{e^{{t^2}}}dt = } \sqrt \pi  \frac{{ - {e^{ - {z^2}}} + w\left( z \right)}}{{2i}},
$$
the Fresnel integral \cite{Abramowitz1972, McKenna1984}
\[
\begin{aligned}
F_r\left( z \right) &= \int\limits_{0}^z  {{e^{i\left( {\pi /2} \right){t^2}}}dt} \\
 &= \left( {1 + i} \right)\left[ {1 - {e^{i\left( {\pi /2} \right){z^2}}}w\left( {\sqrt \pi  \left( {1 + i} \right)z/2} \right)} \right]/2
\end{aligned}
\]
and the normal distribution function \cite{Weisstein2003}
\[	
\begin{aligned}
\Phi \left( z \right) = \frac{1}{{\sqrt {2\pi } }}\int\limits_0^z {{e^{ - {t^2}/2}}dt = \frac{1}{2}{\rm{erf}}\left( {\frac{z}{{\sqrt 2 }}} \right)} \\
 = \frac{1}{2}\left[ {1 - {e^{ - {z^2}/2}}w\left( {\frac{{iz}}{{\sqrt 2 }}} \right)} \right].
\end{aligned}
\]

It is not difficult to show that the complex error function can be represented in an alternative form (see equation (3) in \cite{Srivastava1987} and \cite{Srivastava1992}, see also Appendix A in \cite{Abrarov2015b} for derivation)
\begin{equation}\label{eq_6}
w\left( {x,y} \right) = \frac{1}{{\sqrt \pi  }}\int\limits_0^\infty  {\exp \left( { - {t^2}/4} \right)\exp \left( { - yt} \right)\exp \left( {ixt} \right)dt}.
\end{equation}
This representation of the complex error function will be used for derivation of a rational approximation. 

\subsection{Rational approximation}
In our recent publications we have shown a new technique to obtain a rational approximation for the integrals of kind \cite{Abrarov2014, Abrarov2015c}
$$
\int \limits_0^\infty {{e^{ - {t^2}}}f\left( t \right)dt}.
$$
We apply this approach together with the Fourier transform methodology discussed above in the Introduction.

We can use either of equation \eqref{eq_3a} or \eqref{eq_3b}. For example, we may choose the equation \eqref{eq_3b} corresponding to the inverse Fourier transform. Consider the function $f\left( t \right) = {e^{ - {t^2}/4}}$. Let us find first its forward Fourier transform by substituting $f\left( t \right) = {e^{ - {t^2}/4}}$ into equation \eqref{eq_1a}. These leads to
$$F\left( \nu  \right) = \int\limits_{ - \infty }^\infty  {{e^{ - {t^2}/4}}{e^{ - 2\pi ivt}}dt}  = 2\sqrt \pi  {e^{ - {{\left( {2\pi \nu } \right)}^2}}}.
$$
Now substituting $2\sqrt \pi  {e^{ - {{\left( {2\pi \nu } \right)}^2}}}$ into equation \eqref{eq_3b} yields the following approximation for the exponential function
$$
{e^{ - {t^2}/4}} \approx 2\sqrt \pi  {h_i}\sum\limits_{n =  - N}^N {{e^{ - {{\left( {2\pi n{h_i}} \right)}^2}}}{e^{2\pi itn{h_i}}}},		\qquad - \frac{1}{{2{h_i}}} \le t \le \frac{1}{{2{h_i}}},
$$
or
\begin{equation}\label{eq_7}
{e^{ - {t^2}/4}} \approx 2\sqrt \pi  {h_i}\sum\limits_{n =  - N}^N {{e^{ - {{\left( {2\pi n{h_i}} \right)}^2}}}\cos \left( {2\pi tn{h_i}} \right)},	\qquad - \frac{1}{{2{h_i}}} \le t \le \frac{1}{{2{h_i}}}.
\end{equation}
Taking into account that
$$
{e^{ - {{\left( {2\pi 0{h_i}} \right)}^2}}}\cos \left( {2\pi t0{h_i}} \right) = 1
$$
and
$$
\sum\limits_{n =  - N}^{ - 1} {{e^{ - {{\left( {2\pi n{h_i}} \right)}^2}}}\cos \left( {2\pi tn{h_i}} \right)}  = \sum\limits_{n = 1}^N {{e^{ - {{\left( {2\pi n{h_i}} \right)}^2}}}\cos \left( {2\pi tn{h_i}} \right)}
$$
the approximation \eqref{eq_7} can be simplified as given by
\small
\begin{equation}\label{eq_8}
{e^{ - {t^2}/4}} \approx 2\sqrt \pi  {h_i}\left[ {1 + 2\sum\limits_{n = 1}^N {{e^{ - {{\left( {2\pi n{h_i}} \right)}^2}}}\cos \left( {2\pi tn{h_i}} \right)} } \right],	\qquad - \frac{1}{{2{h_i}}} \le t \le \frac{1}{{2{h_i}}}.
\end{equation}
\normalsize

The right side limitation $t \le 1/\left( {2{h_i}} \right)$ along the positive $t$-axis in equation \eqref{eq_8} can be readily excluded by multiplying both its sides to $\exp \left( { - \sigma t} \right)$ if a constant $\sigma $ is positive and sufficiently large. This can be explained by considering Fig. 1 that shows two functions computed according to right side of equation \eqref{eq_8} at $\sigma  = 0.1$ (blue curve) and $\sigma  = 0.2$ (red curve). For example, at $\sigma  = 0.1$ we can observe two additional peaks at $1/{h_i}$ and $2/{h_i}$ (blue curve). However, as $\sigma $ increases the additional peaks are suppressed stronger to zero due to multiplication to the damping exponential function $\exp \left( { - \sigma t} \right)$. As a result, at $\sigma  = 0.2$ only a single additional peak at $1/{h_i}$ remains visible (red curve). By $\sigma \gtrsim 1$  all additional peaks completely vanish and, therefore, do not contribute to error in integration. Consequently, if the constant $\sigma $ is large enough, say approximately equal or greater than $1$, we can write the approximation
$$
{e^{ - {t^2}/4}}{e^{ - \sigma t}} \approx 2\sqrt \pi  {h_i}\left[ {1 + 2\sum\limits_{n = 1}^N {{e^{ - {{\left( {2\pi n{h_i}} \right)}^2}}}\cos \left( {2\pi tn{h_i}} \right)} } \right]{e^{ - \sigma t}}, \qquad \sigma  \mathbin{\lower.3ex\hbox{$\buildrel>\over
{\smash{\scriptstyle\sim}\vphantom{_x}}$}} 1,
$$
that remains always valid without any limitation along the positive $t$-axis. Assuming $y \ge 0$ we, therefore, can write now
\small
\begin{equation}\label{eq_9}
{e^{ - {t^2}/4}}{e^{ - \left( {y + \sigma } \right)t}} \approx 2\sqrt \pi  {h_i}\left[ {1 + 2\sum\limits_{n = 1}^N {{e^{ - {{\left( {2\pi n{h_i}} \right)}^2}}}\cos \left( {2\pi tn{h_i}} \right)} } \right]{e^{ - \left( {y + \sigma } \right)t}}, \qquad \sigma  \mathbin{\lower.3ex\hbox{$\buildrel>\over
{\smash{\scriptstyle\sim}\vphantom{_x}}$}} 1.
\end{equation}
\normalsize

%\newpage
\begin{figure}[ht]
\begin{center}
\includegraphics[width=24pc]{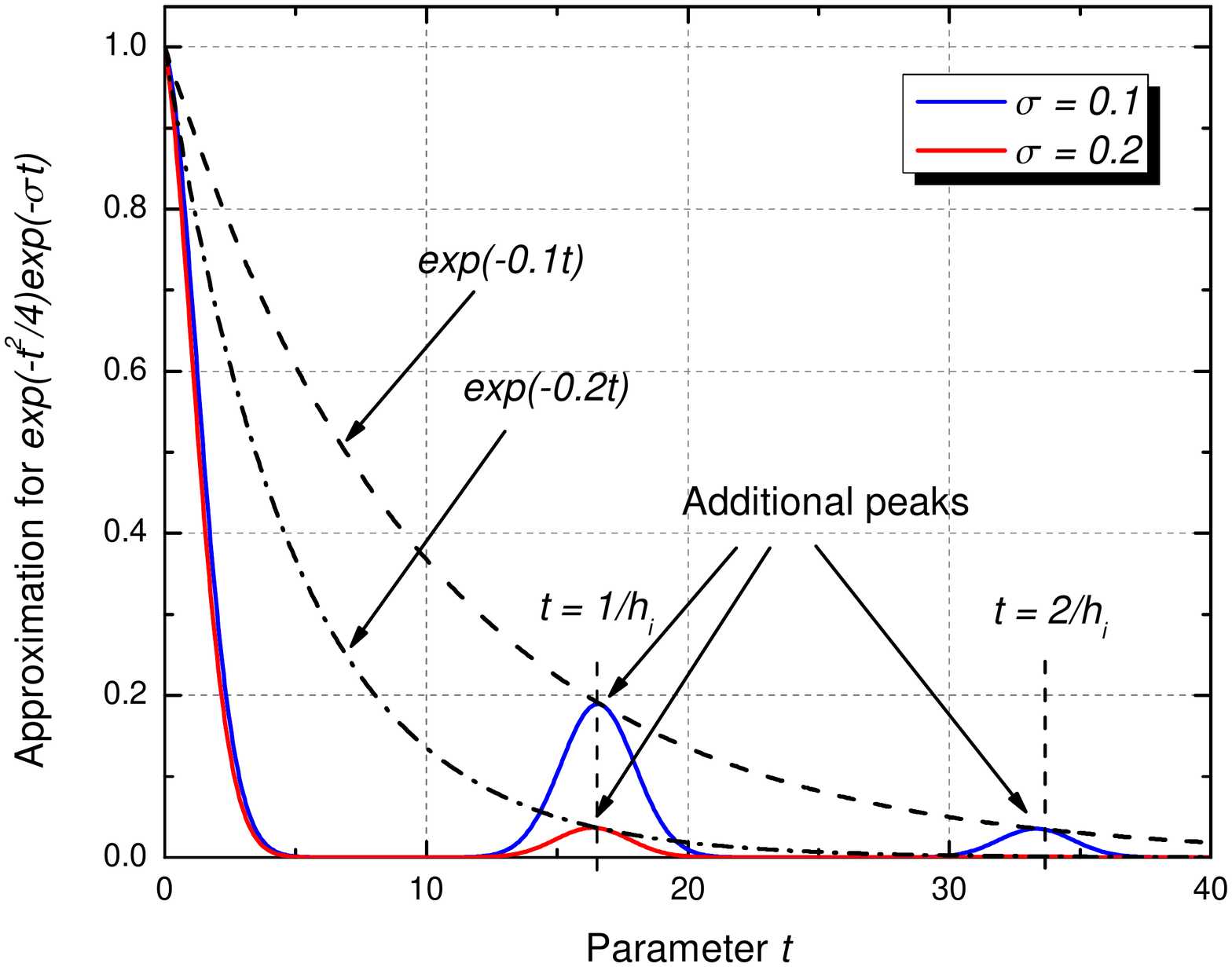}\hspace{1pc}%
\begin{minipage}[b]{28pc}
\vspace{0.3cm}
{\sffamily {\bf{Fig. 1.}} Function approximation for ${e^{ - {t^2}/4}}{e^{ - \sigma t}}$ at $\sigma  = 0.1$ (blue curve) and $\sigma  = 0.2$ (red curve). The dashed and dot-dashed curves are the exponential functions ${e^{ - 0.1t}}$ and ${e^{ - 0.2t}}$, respectively.}
\end{minipage}
\end{center}
\end{figure}

Since ${e^{ - {t^2}/4}}{e^{ - yt}} = {e^{{\sigma ^2}}}{e^{ - {{\left( {t - 2\sigma } \right)}^2}/4}}{e^{ - \left( {y + \sigma } \right)t}}$ from approximation \eqref{eq_9} we obtain
\footnotesize
\begin{equation}\label{eq_10}
{e^{ - {t^2}/4}}{e^{ - yt}} \approx 2\sqrt \pi  {h_i}{e^{{\sigma ^2}}}\left[ {1 + 2\sum\limits_{n = 1}^N {{e^{ - {{\left( {2\pi n{h_i}} \right)}^2}}}\cos \left( {2\pi n{h_i}\left( {t - 2\sigma } \right)} \right)} } \right]{e^{ - \left( {y + \sigma } \right)t}}, \quad \sigma  \mathbin{\lower.3ex\hbox{$\buildrel>\over
{\smash{\scriptstyle\sim}\vphantom{_x}}$}} 1.
\end{equation}
\normalsize
Once again, due to presence of the rapidly damping exponential multiplier ${e^{ - \left( {y + \sigma } \right)t}}$ this approximation is valid without any limitation along the positive $t$-axis. As the peak of the function ${e^{ - {{\left( {t - 2\sigma } \right)}^2}/4}}$ is shifted towards right with respect to the origin, we may regard to the value $\sigma $ as the shift constant.

Finally, substituting approximation \eqref{eq_10} into integral \eqref{eq_6} yields
\begin{equation}\label{eq_11}
w\left( z \right) = i\frac{{2{h_i}{e^{{\sigma ^2}}}}}{z + i\sigma } + \sum\limits_{n = 1}^N {\frac{{{A_n} - i\left( {z + i\sigma } \right){B_n}}}{{C_n^2 - {{\left( {z + i\sigma } \right)}^2}}}},
\end{equation}
where
$$
{A_n} = 8\pi h_i^2n{e^{{\sigma ^2} - {{\left( {2\pi {h_i}n} \right)}^2}}}\sin \left( {4\pi {h_i}n\sigma } \right),
$$
$$
{B_n} = 4{h_i}{e^{{\sigma ^2} - {{\left( {2\pi {h_i}n} \right)}^2}}}\cos \left( {4\pi {h_i}n\sigma } \right)
$$
and 
$$
{C_n} = 2\pi {h_i}n.
$$
As the expansion coefficients ${A_n}$, ${B_n}$ and ${C_n}$ are independent of the argument $z$, the obtained equation \eqref{eq_11} is a rational approximation.

In algorithmic implementation it is more convenient to use $\psi $-function defined as
\begin{equation}
\begin{aligned}\label{eq_12}
\psi \left( z \right) = i\frac{{2{h_i}{e^{{\sigma ^2}}}}}{z} + \sum\limits_{n = 1}^N {\frac{{{A_n} - iz{B_n}}}{{C_n^2 - {z^2}}}} \\
 \Rightarrow w\left( z \right) \approx \psi \left( {z + i\sigma } \right).
\end{aligned}
\end{equation}

\subsection{Computational procedure and error analysis}
Due to a remarkable identity of the complex error function \cite{McKenna1984, Zaghloul2011}
\begin{equation}\label{eq_13}
w\left( { - z} \right) = 2{e^{ - {z^2}}} - w\left( z \right),
\end{equation}
it is sufficient to consider only I and II quadrants in order to cover the entire complex plane. This can be seen explicitly by representation of the identity \eqref{eq_13} in form
$$
w\left( { \pm x, - \left| y \right|} \right) = 2{e^{ - {{\left( { \mp x + i\left| y \right|} \right)}^2}}} - w\left( { \mp x, + \left| y \right|} \right).
$$
Thus, if the parameter $y$ is negative we can simply take it by absolute value and then compute the complex error function according to right side of this equation. Therefore, further we will always assume that $y \ge 0$.

When the argument $z$ is large enough by absolute value, say $\left| {x + iy} \right| \mathbin{\lower.3ex\hbox{$\buildrel>\over
{\smash{\scriptstyle\sim}\vphantom{_x}}$}} 15$, we can truncate the Laplace continued fraction \cite{Gautschi1970, Jones1988}
$$
w\left( z \right) = \frac{{{\mu _0}}}{{z - }}\frac{{1/2}}{{z - }}\frac{1}{{z - }}\frac{{3/2}}{{z - }}\frac{2}{{z - }}\frac{{5/2}}{{z - }}\frac{3}{{z - }}\frac{{7/2}}{{z - }}...\,,	\quad {\mu _0} = i/\pi .
$$
Approximation based on the Laplace continued fraction is rapid in computation. However, its accuracy deteriorates as the argument $z$ decreases by absolute value.

There are different approximations for computation of the narrow-band domain $0 \le x \le 15$ and $0 \le y < {10^{ - 6}}$ \cite{Abrarov2014, Amamou2013, Abrarov2015d}. We can apply, for example, an approximation proposed in our recent work \cite{Abrarov2014}
\small
$$
w\left( {x,y <  < 1} \right) \approx \left( {1 - \frac{y}{{{y_{\min }}}}} \right){e^{ - {x^2}}} + \frac{y}{{{y_{\min }}}}K\left( {x,{y_{\min }}} \right) + iL\left( {x,{y_{\min }}} \right), \qquad y_{\min} <  < 1,
$$
\normalsize
where $L \left( x, y_{\min} \right) = {\mathop{\rm Im}\nolimits} \left[ w \left( x, y_{\min} \right) \right]$ and ${y_{\min }}$ can be taken equal to ${10^{ - 5}}$. It has been shown that this approximation can provide accuracy better than ${10^{ - 9}}$ over the narrow-band domain $0 \le x \le 15$ and $0 \le y < {10^{ - 6}}$.

The domain $\left| {x + iy} \right| \le 15 \cap y \ge {10^{ - 6}}$ is the most difficult for computation. Nevertheless, with only $17$ summation terms (at $N = 16$) the proposed rational approximation \eqref{eq_12} covers this domain providing high-accuracy and rapid computation. In computational procedure we have to choose properly the margin value ${\nu _m}$ for the exponential function ${e^{ - {{\left( {2\pi \nu } \right)}^2}}}$ that appears from the forward Fourier transform $F\left( \nu  \right) = 2\sqrt \pi  {e^{ - {{\left( {2\pi \nu } \right)}^2}}}$. As it has been justified by Melone {\it{et al.}} \cite{Milone1988}, the margin value for integration involving the exponential function ${e^{ - {t^2}}}$ can be taken as $t = {t_m} = 6$. We can use this result in order to determine the required value by solving the following equation with respect to the variable $\nu $ as follows
$$
{e^{ - {{\left( {2\pi \nu } \right)}^2}}} = {\left. {{e^{ - {t^2}}}} \right|_{t = 6}} \Rightarrow {\left( {2\pi \nu } \right)^2} = 36.
$$
There are two solutions for this equation ${\nu _{1,2}} =  \pm 6/\left( {2\pi } \right)$. Consequently, the margin value for the exponential function ${e^{ - {{\left( {2\pi \nu } \right)}^2}}}$ can be taken as ${\nu _m} = 6/\left( {2\pi } \right)$. As a parameter ${h_i}$ is the step between two adjacent sampling points along positive $\nu $-axis (see \cite{Abrarov2015a} for details), its value can be calculated as ${h_i} = {\nu _m}/N$. Taking $N = 16$ we can find that ${h_i} = {\nu _m}/16 \approx {\rm{5}}{\rm{.968310365946075}} \times {\rm{1}}{{\rm{0}}^{ - 2}}$.

In order to quantify the accuracy of the rational approximation \eqref{eq_12} we may define the relative errors
$$
{\Delta _{{\mathop{\rm Re}\nolimits} }} = \left| {\frac{{{\mathop{\rm Re}\nolimits} \left[ {{w_{ref}}\left( {x,y} \right)} \right] - {\mathop{\rm Re}\nolimits} \left[ {w\left( {x,y} \right)} \right]}}{{{\mathop{\rm Re}\nolimits} \left[ {{w_{ref}}\left( {x,y} \right)} \right]}}} \right|
$$
and
$$
{\Delta _{{\mathop{\rm Im}\nolimits} }} = \left| {\frac{{{\mathop{\rm Im}\nolimits} \left[ {{w_{ref}}\left( {x,y} \right)} \right] - {\mathop{\rm Im}\nolimits} \left[ {w\left( {x,y} \right)} \right]}}{{{\mathop{\rm Im}\nolimits} \left[ {{w_{ref}}\left( {x,y} \right)} \right]}}} \right|,
$$
where ${w_{ref}}\left( {x,y} \right)$ is the reference, for the real and imaginary parts, respectively. The highly accurate reference values can be generated by using, for example, Algorithm 680 \cite{Poppe1990a, Poppe1990b}, recently published Algorithm 916 \cite{Zaghloul2011} or C++ code from the RooFit package, CERN\text{'}s library \cite{Karbach2014}.

Figure 2 shows ${\log _{10}}{\Delta _{{\mathop{\rm Re}\nolimits} }}$ for the real part of the complex error function computed over the domain $0 \le x \le 15$ and ${10^{ - 6}} \le y \le 15$ at $N = 16$, $\sigma  = 1.5$  and ${h_i} = {\rm{5}}{\rm{.968310365946075}} \times {\rm{1}}{{\rm{0}}^{ - 2}}$. As we can see from this figure, the rational approximation \eqref{eq_12} provides accuracy ${10^{ - 15}}$ (blue color) over the most of this domain. Although accuracy deteriorates with decreasing $y$, it remains better than ${10^{ - 9}}$ (red color) in the range ${10^{ - 4}} \le y \le {10^{ - 6}}$. This indicates that at the same $N = 16$ the accuracy of the rational approximation \eqref{eq_12} is by several orders of the magnitude higher than the accuracy of the Weideman\text{'}s rational approximation (see equation 38(I) in \cite{Weideman1994}).

%\newpage
\begin{figure}[ht]
\begin{center}
\includegraphics[width=22pc]{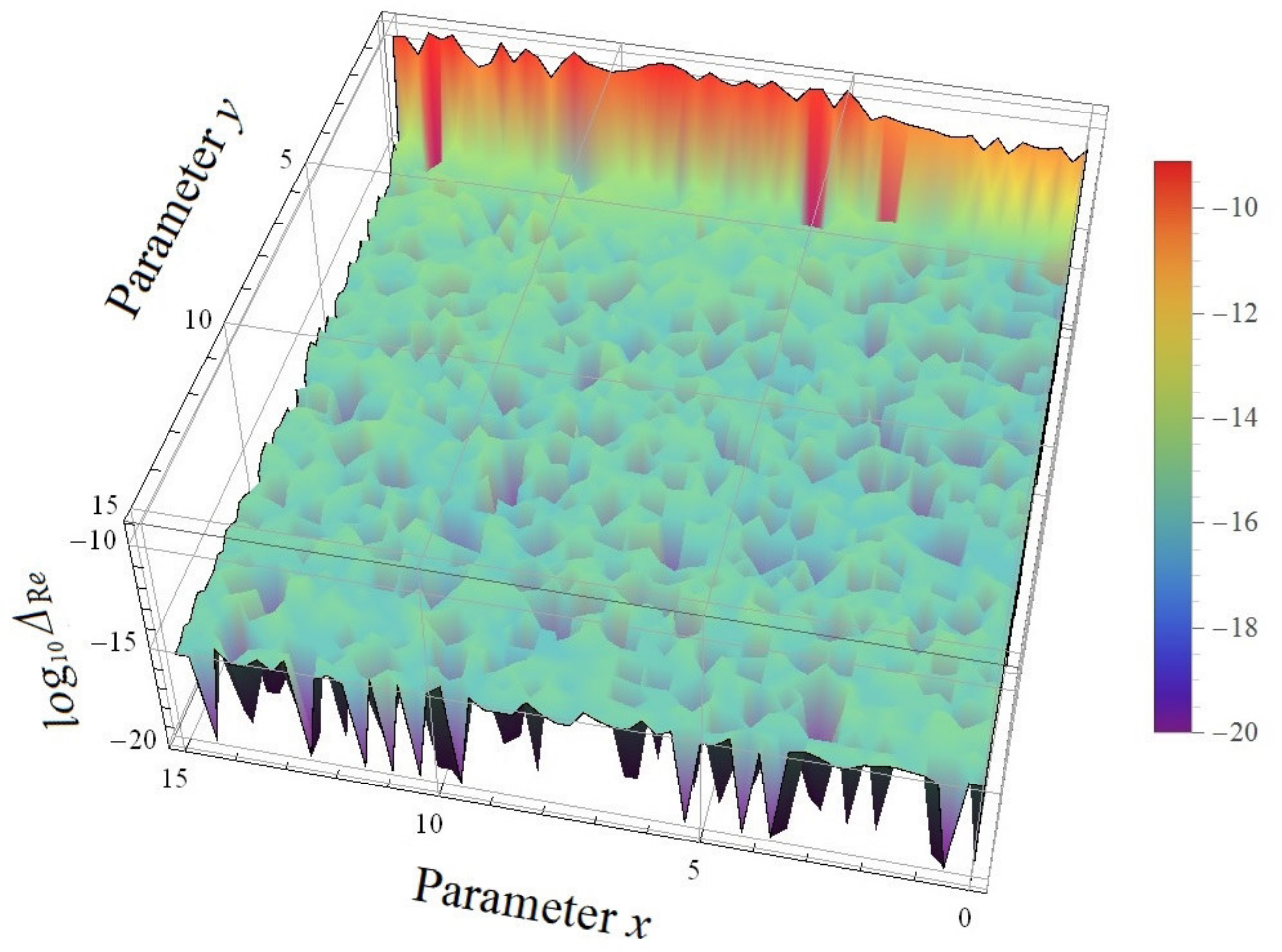}\hspace{2pc}%
\begin{minipage}[b]{28pc}
\vspace{0.3cm}
{\sffamily {\bf{Fig. 2.}} The logarithm of the relative error ${\log _{10}}{\Delta _{{\mathop{\rm Re}\nolimits} }}$ for the real part of the rational approximation \eqref{eq_12} over the domain $0 \le x \le 15 \cap {10^{ - 6}} \le y \le {10^4}$.}
\end{minipage}
\end{center}
\end{figure}

Figure 3 illustrates ${\log _{10}}{\Delta _{{\mathop{\rm Im}\nolimits} }}$ for the imaginary part of the complex error function also computed over the domain $0 \le x \le 15$ and ${10^{ - 6}} \le y \le 15$ at $N = 16$, $\sigma  = 1.5$ and ${h_i} = {\rm{5}}{\rm{.968310365946075}} \times {\rm{1}}{{\rm{0}}^{ - 2}}$. One can see that in the imaginary part the accuracy is also highly accurate ${10^{ - 15}}$ (blue color) over the most domain. There is only a small area $0 \le x < 1$ and ${10^{ - 6}} \le y \le {10^{ - 4}}$ near the origin where the accuracy deteriorates as the parameters $x$ and $y$ both tend to zero. Nevertheless, the accuracy in this area still remains high and better than ${10^{ - 9}}$ (red color).

%\newpage
\begin{figure}[ht]
\begin{center}
\includegraphics[width=22pc]{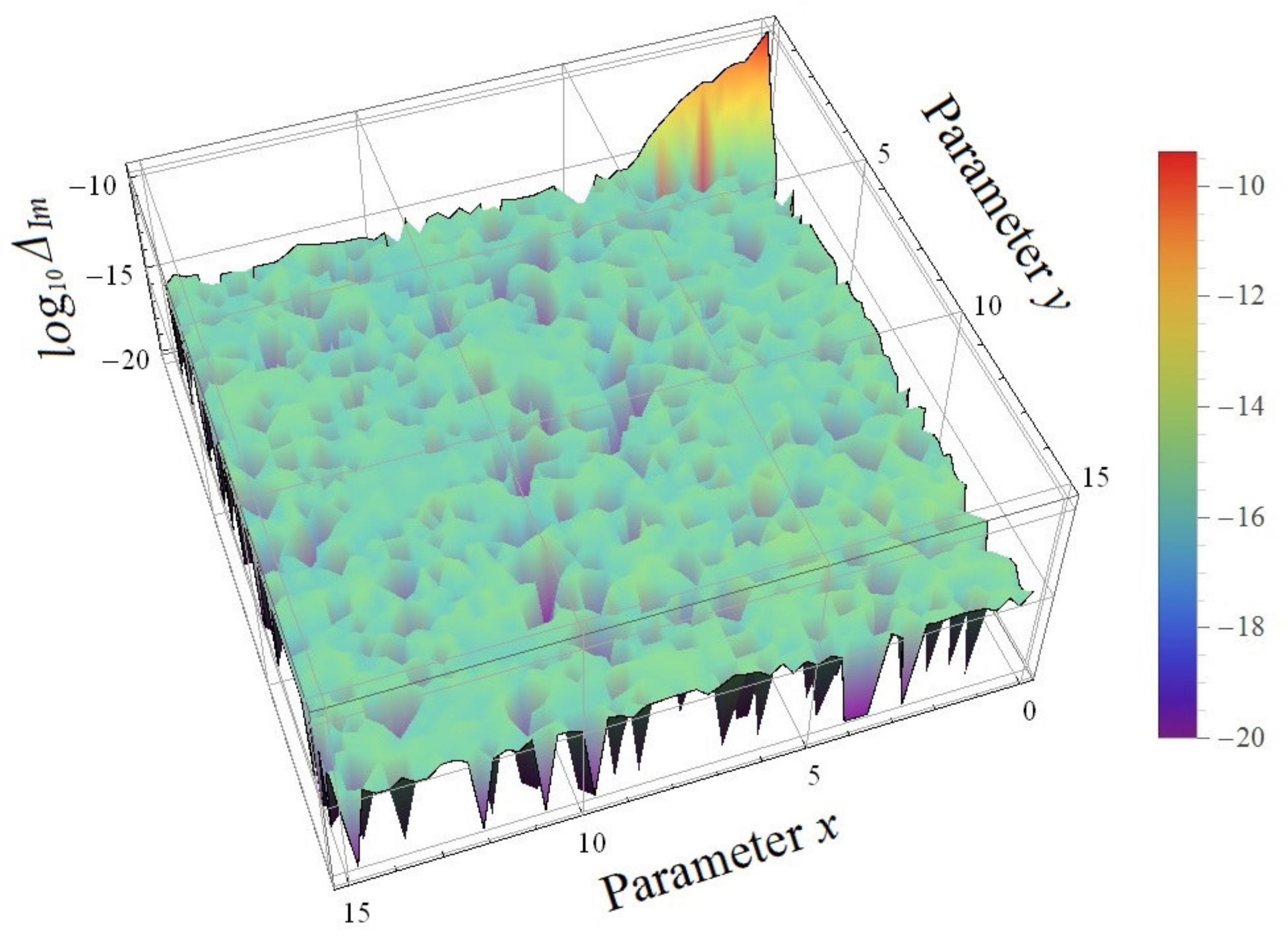}\hspace{2pc}%
\begin{minipage}[b]{28pc}
\vspace{0.3cm}
{\sffamily {\bf{Fig. 3.}} The logarithm of the relative error ${\log _{10}}{\Delta _{{\mathop{\rm Im}\nolimits} }}$ for the imaginary part of the rational approximation \eqref{eq_12} over the domain $0 \le x \le 15 \cap {10^{ - 6}} \le y \le {10^4}$.}
\end{minipage}
\end{center}
\end{figure}

The computational test reveals that with only $17$ summation terms (at $N = 16$) the rational approximation \eqref{eq_12} alone can cover the entire HITRAN domain $0 \le x \le 40,000 \cap {10^{ - 4}} \le y \le {10^2}$ providing average accuracy ${10^{ - 15}}$ for an input array consisting of $3 \times 10^7$ elements. Algorithmic implementation of the rational approximation \eqref{eq_12} results to the same computational speed as that of described in our recent work where we proposed a sampling by incomplete cosine expansion of the sinc function to approximate the complex error function \cite{Abrarov2015c}. 

A Matlab subroutine code that covers the HITRAN domain with high-accuracy is presented in Appendix A.
	
\section{Conclusion}
We present a new efficient rational approximation to the complex error function by application of the Fourier transform that provides computationally rapid and highly accurate results. The computational test we performed with only $17$ summation terms shows that the accuracy of the rational approximation of the complex error function is ${10^{ - 15}}$ over the most domain of practical importance. In particular, the proposed rational approximation of the complex error function alone can cover with high accuracy the entire domain $0 \le x \le 40,000 \cap {10^{ - 4}} \le y \le {10^2}$ required for the HITRAN-based spectroscopic applications.

\section*{Acknowledgments}
This work is supported by National Research Council Canada, Thoth Technology Inc. and York University.

\section*{Appendix A}
\footnotesize
\begin{verbatim}
function VF = comperf(z)
 
% This subroutine function file computes the complex error function, also 
% known as the Faddeeva function. It covers the entire HITRAN domain 
% 0 <= x <= 40,000 and 10^-4 <= y <= 10^2. However, it may be used only in 
% the most difficult domain |x + 1i*y| <= 15 and y > = 10^-6. See the
% article that describes how the entire complex plain can be covered.
 
% The code is written by Sanjar M. Abrarov and Brendan M. Quine, York 
% University, October, 2015.
 
if any(imag(z) < 10^-6)
    disp('One or more imag(z) is less than 10^-6. Computation terminated.')
    VF = NaN;
    return
end
 
num = 16; % number of summation terms is 16 + 1 = 17
vm = 6/(2*pi); % margin value
hi = vm/num; % sampling step
sig = 1.5; % the shift constant
 
n = 1:num; % define array n

% Define the expansion coefficients
An = 8*pi*hi^2*n.*exp(sig^2 - (2*pi*hi*n).^2).*sin(4*pi*hi*n*sig);
Bn = 4*hi*exp(sig^2 - (2*pi*hi*n).^2).*cos(4*pi*hi*n*sig);
Cn = 2*pi*hi*n;

z = z + 1i*sig; % redefine input z (see formula (12) for the psi-function)
zz = z.^2;
 
    VF = 1i*(2*hi*exp(sig^2))./z; % define first term
    for n = 1:num
        VF = VF + (An(n) - 1i*z*Bn(n))./(Cn(n)^2 - zz);
    end
end
\end{verbatim}
\normalsize

%\bigskip
%\newpage


\begin{thebibliography}{9}

\bibitem{Bracewell2000}
R.N. Bracewell, The Fourier transform and its application, $3^\text{rd}$ ed., McGraw-Hill, 2000.

\bibitem{Hansen2014}
E.W. Hansen, Fourier transforms. Principles and applications, John Wiley \& Sons, 2014.

\bibitem{Mishra2012}
V.N. Mishra and L.N. Mishra, Trigonometric approximation of signals (functions) in $L_{p}$-norm. Int. J. Contemp. Math. Sci., 7 (19) (2012) 909-918. \url{http://www.m-hikari.com/ijcms/ijcms-2012/17-20-2012/narayanmishraIJCMS17-20-2012.pdf}

\bibitem{Mishra2013}
V.N. Mishra, H.H. Khan, I.A. Khan, K. Khatri and L.N. Mishra, Trigonometric approximation of signals (functions) belonging to the $Lip \left( \xi \left( t \right),r \right)$, $\left( r > 1 \right)$-class by $\left( E,q \right)$ $\left( q > 0 \right)$-means of the conjugate series of its Fourier series. Advan. Pure Math., 3 (2013) 353-358. \url{http://dx.doi.org/10.4236/apm.2013.33050}

\bibitem{Boyd2013}
J.P. Boyd, A comparison of companion matrix methods to find roots of a trigonometric polynomial, J. Comp. Phys., 246 (2013) 96–112. \url{http://dx.doi.org/10.1016/j.jcp.2013.03.022}

\bibitem{Abrarov2015a}
S.M. Abrarov and B. M. Quine, Representation of the Fourier transform as a weighted sum of the complex error functions, arXiv:1507.01241v3. \url{http://arxiv.org/pdf/1507.01241v3.pdf}

\bibitem{Rothman2013}
L.S. Rothman, I.E. Gordon, Y. Babikov, A. Barbe, D.C. Benner, P.F. Bernath, M. Birk, L. Bizzocchi, V. Boudon, L.R. Brown, A. Campargue, K. Chance, E.A. Cohen, L.H. Coudert, V.M. Devi, B.J. Drouin, A. Fayt, J.-M. Flaud, R.R. Gamache, J.J. Harrison, J.-M. Hartmann, C. Hill, J.T. Hodges, D. Jacquemart, A. Jolly, J. Lamouroux, R.J. Le Roy, G. Li, D.A. Long, O.M. Lyulin, C.J. Mackie, S.T. Massie, S. Mikhailenko, H.S.P. M\"{u}ler, O.V. Naumenko, A.V. Nikitin, J. Orphal, V. Perevalov, A. Perrin, E.R. Polovtseva and C. Richard, The HITRAN2012 molecular spectroscopic database, J. Quant. Spectrosc. Radiat. Transfer, 130 (2013) 4-50. \url{http://dx.doi.org/10.1016/j.jqsrt.2013.07.002}

\bibitem{Faddeyeva1961}
V.N. Faddeyeva, and N.M. Terent\text{'}ev, Tables of the probability integral $w\left( z \right) = {e^{ - {z^2}}}\left( {1 + \frac{{2i}}{{\sqrt \pi  }}\int_0^z {{e^{{t^2}}}dt} } \right)$ for complex argument. Pergamon Press, Oxford, 1961.

\bibitem{Gautschi1970}
W. Gautschi, Efficient computation of the complex error function. SIAM J. Numer. Anal., 7 (1970) 187-198. \url{http://dx.doi.org/10.1137/0707012}

\bibitem{Abramowitz1972}
M. Abramowitz and I.A. Stegun. Error Function and Fresnel Integrals. Handbook of mathematical functions with formulas, graphs, and mathematical tables. $9^{\text{th}}$ ed. New York 1972, 297-309.

\bibitem{Armstrong1972}
B.H. Armstrong and B.W. Nicholls, Emission, absorption and transfer of radiation in heated atmospheres. Pergamon Press, New York, 1972.

\bibitem{Schreier1992}
F. Schreier, The Voigt and complex error function: A comparison of computational methods. J. Quant. Spectrosc. Radiat. Transfer, 48 (1992) 743-762. \url{http://dx.doi.org/10.1016/0022-4073(92)90139-U}

\bibitem{Weideman1994}
J.A.C. Weideman, Computation of the complex error function. SIAM J. Numer. Anal., 31 (1994) 1497-1518. \url{http://dx.doi.org/10.1137/0731077}

\bibitem{Letchworth2007}
K.L. Letchworth and D.C. Benner, Rapid and accurate calculation of the Voigt function, J. Quant. Spectrosc. Radiat. Transfer, 107 (2007) 173-192. \url{http://dx.doi.org/10.1016/j.jqsrt.2007.01.052}

\bibitem{Pagnini2010}
G. Pagnini and F. Mainardi, Evolution equations for the probabilistic generalization of the Voigt profile function, J. Comput. Appl. Math., 233 (2010) 1590-1595. \url{http://dx.doi.org/10.1016/j.cam.2008.04.040}

\bibitem{Abrarov2010}
S.M. Abrarov, B.M. Quine and R.K. Jagpal, High-accuracy approximation of the complex probability function by Fourier expansion of exponential multiplier,  Comp. Phys. Commun., 181 (5) (2010) 876-882. \url{http://dx.doi.org/10.1016/j.cpc.2009.12.024}

\bibitem{Edwards1992}
D.P. Edwards, GENLN2: A general line-by-line atmospheric transmittance and radiance model, NCAR technical note, 1992. \url{http://dx.doi.org/10.5065/D6W37T86}

\bibitem{Quine2002}
B.M. Quine and J.R. Drummond, GENSPECT: a line-by-line code with selectable interpolation error tolerance J. Quant. Spectrosc. Radiat. Transfer 74 (2002) 147-165. \url{http://dx.doi.org/10.1016/S0022-4073(01)00193-5}

\bibitem{Christensen2012}
L.E. Christensen, G.D. Spiers, R.T. Menzies and J.C Jacob, Tunable laser spectroscopy of $\rm{CO}_{2}$ near $2.05 \, {\mu}m$: Atmospheric retrieval biases due to neglecting line-mixing, J. Quant. Spectrosc. Radiat. Transfer, 113 (2012) 739-748. \url{http://dx.doi.org/10.1016/j.jqsrt.2012.02.031}

\bibitem{Berk2013}
A. Berk, Voigt equivalent widths and spectral-bin single-line transmittances: Exact expansions and the MODTRAN{\circledR}5 implementation, J. Quant. Spectrosc. Radiat. Transfer, 118 (2013) 102-120. \url{http://dx.doi.org/10.1016/j.jqsrt.2012.11.026}

\bibitem{Quine2013}
B.M. Quine and S.M. Abrarov, Application of the spectrally integrated Voigt function to line-by-line radiative transfer modelling. J. Quant. Spectrosc. Radiat. Transfer, 127 (2013) 37-48. \url{http://dx.doi.org/10.1016/j.jqsrt.2013.04.020}

\bibitem{Fried1961}
B.D. Fried and S.D. Conte. The plasma dispersion function. New York: Academic Press, 1961.

\bibitem{Cody1970}
W.J. Cody, K.A. Paciorek and H.C. Thacher, Chebyshev approximations for Dawson\text{'}s integral. Math. Comp. 24 (1970) 171-178. \url{http://dx.doi.org/10.1090/S0025-5718-1970-0258236-8}

\bibitem{McCabe1974}
J.H. McCabe, A continued fraction expansion with a truncation error estimate for Dawson\text{'}s integral, Math. Comp. 28 (1974) 811-816. \url{http://dx.doi.org/10.1090/S0025-5718-1974-0371020-3}

\bibitem{Rybicki1989}
G.B. Rybicki, Dawson\text{'}s integral and the sampling theorem, Comp. Phys., 3 (1989) 85-87. \url{http://dx.doi.org/10.1063/1.4822832}

\bibitem{Boyd2008}
J.P. Boyd, Evaluating of Dawson\text{'}s integral by solving its differential equation using orthogonal rational Chebyshev functions, Appl. Math. Comput., 204 (2) (2008) 914-919. \url{http://dx.doi.org/10.1016/j.amc.2008.07.039}

\bibitem{Abrarov2015b}
S.M. Abrarov and B.M. Quine, A rational approximation for the Dawson\text{'}s integral of real argument, arXiv:1505.04683. \url{http://arxiv.org/pdf/1505.04683.pdf}

\bibitem{McKenna1984}
S.J. McKenna, A method of computing the complex probability function and other related functions over the whole complex plane. Astrophysics and Space Science, 107 (1) (1984) 71-83. \url{http://dx.doi.org/10.1007/BF00649615}

\bibitem{Weisstein2003}
E.W. Weisstein, CRC concise encyclopedia of mathematics. Chapman \& Hall/CRC, $2^\text{{nd}}$ ed. 2003.

\bibitem{Srivastava1987}
H.M. Srivastava and E.A. Miller, A unified presentations of the Voigt functions, Astrophys. Space Sci., 135 (1987) 111-118. \url{http://dx.doi.org/10.1007/BF00644466}

\bibitem{Srivastava1992}
H.M. Srivastava and M.P. Chen, Some unified presentations of the Voigt functions. Astrophys. Space Sci., 192 (1) (1992) 63-74. \url{http://dx.doi.org/10.1007/BF00653260}

\bibitem{Abrarov2014}
S.M. Abrarov and B.M. Quine, Master-slave algorithm for highly accurate and rapid computation of the Voigt/complex error function, J. Math. Research, 6 (2) (2014) 104-119. \url{http://dx.doi.org/10.5539/jmr.v6n2p104}

\bibitem{Abrarov2015c}
S.M. Abrarov and B.M. Quine, Sampling by incomplete cosine expansion of the sinc function: Application to the Voigt/complex error function, Appl. Math. Comput., 258 (2015) 425-435. \url{http://dx.doi.org/ 10.1016/j.amc.2015.01.072}

\bibitem{Zaghloul2011}
M.R. Zaghloul and A.N. Ali, Algorithm 916: computing the Faddeyeva and Voigt functions. ACM Transactions on Mathematical Software, 38 (2011) 15:1-15:22. \url{http://dx.doi.org/10.1145/2049673.2049679}

\bibitem{Jones1988}
W.B. Jones and W.J. Thron, Continued fractions in numerical analysis. Appl. Num. Math., 4(2-4) (1988) 143-230. \url{http://dx.doi.org/10.1016/0168-9274(83)90002-8}

\bibitem{Amamou2013}
H. Amamou, B. Ferhat and A. Bois, Calculation of the Voigt Function in the region of very small values of the parameter $a$ where the calculation is notoriously difficult, Amer. J. Anal. Chem., 4 (2013) 725-731. \url{http://dx.doi.org/10.4236/ajac.2013.412087}

\bibitem{Abrarov2015d}
S.M. Abrarov and S.M. Quine, Accurate approximations for the complex error function with small imaginary argument, J. Math. Research 7 (1) (2015) 44-53. \url{http://dx.doi.org/10.5539/jmr.v1n1p44}

\bibitem{Milone1988}
A.A.E. Milone, L.A. Milone and G.E. Bobato, Numerical evaluation of the line broadening function  . Astrophysics. and Space Science, 147 (2) (1988) 229-234. \url{http://dx.doi.org/10.1007/BF00645667}

\bibitem{Poppe1990a}
G.P.M. Poppe and C.M.J. Wijers, More efficient computation of the complex error function. ACM Transact. Math. Software, 16 (1990) 38-46. \url{http://dx.doi.org/10.1145/77626.77629}

\bibitem{Poppe1990b}
G.P.M. Poppe and C.M.J. Wijers, Algorithm 680: evaluation of the complex error function. ACM Transact. Math. Software, 16 (1990) 47. \url{http://dx.doi.org/10.1145/77626.77630}

\bibitem{Karbach2014}
T.M. Karbach, G. Raven and M. Schiller, Decay time integrals in neutral meson mixing and their efficient evaluation, arXiv:1407.0748. \url{http://arxiv.org/pdf/1407.0748v1.pdf}

\end{thebibliography}
\end{document}